\documentclass[preprint]{elsarticle}
\usepackage{amsmath,amssymb}
\usepackage{pstricks,pst-node}
\begin{document}

\newtheorem{lem}{Lemma}[section]
\newtheorem{thm}{Theorem}[section]
\newproof{pf}{Proof}
\newdefinition{dfn}{Definition}[section]
\title{Automorphisms of the bipartite graph planar algebra}

\author[rdb]{R. D.~Burstein\corref{cor1}}
\cortext[cor1]{Corresponding author.  Tel. 1-615-322-6672
Fax 1-615-343-0215}
\ead{richard.d.burstein@vanderbilt.edu}
\address[rdb]{Department of Mathematics\\
1326 Stevenson Center\\
Vanderbilt University\\
Nashville, TN 37240}

\maketitle

\begin{abstract}
For any abstract subfactor planar algebra $P$, there exists a finite index extremal subfactor $M_0 \subset M_1$
with $P$ as its standard invariant.  In this paper, we classify the automorphism group of a bipartite
graph planar algebra, and obtain subfactor planar subalgebras by taking fixed points under groups of
automorphisms.  This construction provides both new examples of subfactors and new descriptions of
the planar algebras of previously known examples.


\end{abstract}

\section{Introduction}

Planar algebras were introduced by Jones in \cite{JonP}.  They are a powerful tool for studying subfactors,
providing a graphical calculus on the standard invariant of an finite index extremal subfactor of type 
$II_1$ (\cite{JonP},\cite{Pop}).
The planar operad is the set of planar tangles with zero or more internal disks and a checkerboard shading,
with a distinguished region of the boundary and each internal disk, all taken up to isotopy.  
The operation is gluing:
a tangle may be pasted into an internal disk of another tangle (matching up the distinguished boundary
regions) if the number of strands and the shading
are compatible.
A planar algebra is a graded vector space $P = (V_n^{\pm})$, $n \ge 0$, 
along with an associative action of the planar operad. 

If $M_0 \subset M_1$ is a finite index extremal $II_1$ subfactor with Jones tower $M_0 \subset M_1 \subset M_2 \subset ...$,
then we may take $V_n^+ = M_0' \cap M_n$, $V_n^- = M_1' \cap M_{n+1}$.  There is an operad action defined in \cite{JonP}
making this into a planar algebra.  Conversely, Jones describes a certain list of additional properties which make $P$ into
a Popa system, implying that there exists a subfactor $M_0 \subset M_1$ of which $P$ is the standard invariant ~\cite{Pop}.  
A planar algebra with this list of properties is called a subfactor planar algebra (SPA).  Other methods of constructing a subfactor
from an SPA have since been obtained (\cite{GJS},\cite{JSW},\cite{KS}), which use more diagrammatic notation
than Popa's original method.

Constructing an SPA abstractly therefore implies the existence of a corresponding subfactor.  
This method may be used
to find new subfactors (as in \cite{BJ1} or \cite{BMPS}),
or provide new proofs of the existence of subfactors with specified properties.  
Abstract constructions of SPAs of previously known subfactors can provide new insight into the structure of their
standard invariants (e.g. ~\cite{BDG1},\cite{BDG2},\cite{Pet}).

A planar algebra may be constructed from any finite bipartite graph \cite{JonB}, and some infinite graphs as well
\cite{GJS}.  These bipartite graph planar algebras
(BGPAs) are almost never of subfactor type, because their vector spaces are too large.
However, they possess several of the necessary properties required for SPAs, which 
are inherited by planar subalgebras.  We may therefore try to find SPAs by looking at small planar subalgebras
of BGPAs.    

It is generally difficult to show that a graded subspace of a BGPA
is closed under the action of the planar operad, although progress has been made in the single generator case 
(e.g. ~\cite{BJ2},\cite{BJ3},\cite{Pet},\cite{BMPS}).  
In this paper, we describe a method for doing so by taking fixed points of a BGPA
under a group of automorphisms.
In section 2 we describe BGPAs, with particular attention to infinite graphs.
We introduce a slightly different notation from \cite{GJS} and \cite{JonB}; this simplifies the computations of 
section 3, where we compute the automorphism group of an arbitrary BGPA $P_\Gamma$.  In section 4 we find conditions for
a planar subalgebra of a BGPA to be an SPA, and we conclude in section 5 by presenting some examples of SPAs 
obtained by this planar fixed point construction.

\section{The bipartite graph planar algebra}

The bipartite graph planar algebra (BGPA) is described by Jones in \cite{JonB}.  The data required
are a finite bipartite graph $\Gamma$ and a spin vector: i.e. a function from the vertices
of the graph to the positive real numbers.  Given such data,
there is a planar algebra $P_{\Gamma} = (V_n^{\pm}), n \ge 0$,  where each $V_n^{\pm}$ has a basis labelled by 
loops of length $2n$ in the graph, starting at even or odd vertices depending on sign.

We will consider BGPAs on infinite graphs as well.  In this case, the vector spaces
will be infinite dimensional, so there are some new topological considerations.  To simplify
later computations, we will describes
the vector spaces of the planar algebra as operators on a certain Hilbert space. 

In our diagrams, we will use the convention of \cite{JSW} that a thick line represents as
many parallel strands as necessary.  Also we will omit shading in a diagram
when both shadings can occur.  When not otherwise specified, the distinguished boundary region
of a tangle is on the left.

The discussion that follows is taken from \cite{JonB}, adjusted where necessary to allow for
infinite graphs.

Let $\Gamma$ be a locally finite bipartite graph, i.e. each vertex of $\Gamma$ is an endpoint of at most $A$ edges.
Let $\mu$ be a function from the vertices of $\Gamma$
to the positive real numbers obeying the following \textit{local boundedness} condition: there is some $M>0$
such that for any two adjacent vertices $v$ and $w$, we have $\mu(v)/\mu(w)<M$.
Let $l_n^{\pm}$ be the set of loops of length $2n$ on $\Gamma$,
starting at an even $(+)$ or odd $(-)$ vertex.  Then the vector space $V_n^{\pm}$ of the
bipartite graph planar algebra $P_{\Gamma}$ is the set of bounded
functions from $l_n^{\pm}$ to $\mathbb{C}$.

The \textit{boundary type} of a tangle is the ordered pair $(n, \pm)$, where $2n$ strands
intersect the tangle boundary transversely and the distinguished region of the tangle is shaded $(+)$ or
unshaded $(-)$.  The same definition is used for the boundary type of each internal disk of a tangle.

Let $T$ be a planar tangle with $k$ internal disks.  Let the boundary type of $T$ be $(n, \pm)$, and
let the internal disks of $T$ have boundary types respectively $(n_1, \pm), (n_2, \pm), ..., (n_k, \pm)$.  Then to describe
the planar operad, for each set of inputs $(x_1, ..., x_k)$ (with $x_k \in V_{n_k}^{\pm}$) we must assign
the output $Z(T) \in V_n^{\pm}$.  This assignment must agree with gluing, be multilinear in the inputs, and be
isotopy invariant.  We define $Z(T)$ following \cite{JonB}.

A \textit{state} of the tangle is a function 
$\sigma$ which maps the strands of $T$ to edges of $\Gamma$, and the regions of $T$ to vertices of $\Gamma$.
Shaded regions are mapped to positive vertices, and unshaded regions to negative ones.
A state must obey a compatibility condition: if a strand $S$ is adjacent to a region $R$, then $\sigma(R)$
must be one of the endpoints of $\sigma(S)$.

A state $\sigma$ is \textit{compatible} with a given loop if treading the output of $\sigma$ counterclockwise
around the boundary of $T$, starting from the distinguished region, produces that loop.

A \textit{singularity} of $T$ is a local maximum or local minimum of a strand.  Fix a state $\sigma$;
let $v$ be the vertex associated to the concave side of the singularity by $\sigma$, and $w$ the vertex associated
to the convex side.  Then the value of the singularity is $\mu(v)/\mu(w)$.

A state $\sigma$ associates a loop $L_i$ to each internal disk $D_i$ of $T$, obtained by reading the output of $\sigma$
counterclockwise starting from the distinguished region of the disk.  The value of the state on $D_i$ is then 
the value of the $i$th input $x_i$ at this loop $L_i$.

Now we can define the value of $Z(T)$ on a specified loop $L \in l_n^{\pm}$.  This is

$$\sum_{\sigma \, compatible \, with \, L} \left(
\prod_{s \in singularities \, of \, T}\sigma(s) 
\prod_{D_i \in internal \, disks}\sigma(D_i)
\right)
$$

where $\sigma$ is evaluated on disks and internal singularities as above.

By local finiteness of $\Gamma$, only finitely many states are compatible with $L$.  Specifically, let $m$ be the largest
number of strands that must be crossed to get from the distinguished boundary region of $T$ to any other region (boundary or internal).
Then if $v$ is the first vertex of $L$, every compatible state must assign all regions of $T$ to vertices within a distance
of $m$ from $v$ on the graph.  A state is determined by its value on strands.  Let each vertex of $\Gamma$ contact at most $e_{max}$ edges.  
Then each strand must be chosen from at most $e_{max}^m$ possibilities, and if $T$ contains $a$ distinct strands,
each loop is compatible with at most $e_{max}^{am}$ states.

The local boundedness condition on $\mu$ ($\mu(v)/\mu(w) < M < \infty$ for all adjacent vertices $v$ and $w$) means
that each singularity has value at most $M$, for any state. Let $T$ have $b$ singularities.
 Elements of $V_{n_i}^{\pm}$ are bounded; let 
$N$ be the largest bound of any input $x_i$, $1 \le i \le k$.

Putting this together, we find that the tangle output $Z(T)$ may be evaluated on each loop $L$ as a finite sum.  Moreover, this sum is bounded
by $e_{max}^{am} M^b N^k$, so the output is a bounded function on $l_n^{\pm}$ 
and this evaluation rule produces an element of $V_n^{\pm}$.

The proofs of \cite{JonB} that this map is multilinear, isotopy invariant and respects gluing may be used without alteration, 
since local finiteness of $\Gamma$ implies that all necessary sums are finite.  Therefore the above definition of $Z(T)$
produces a planar algebra.

There is a natural antilinear involution, which we refer to as $^*$, of each vector space $V_n^{\pm}$ (see \cite{JonB}).
For each loop $L \in l_n^{\pm}$, the reversed loop $L'$ consists of the same list of vertices and edges,
taken in the opposite order.  
Then the involution is defined by $A^*(L) = \overline{A(L')}$, for all $A \in V_n^{\pm}$ and $L \in l_n^{\pm}$.

Let $\rho_T$ be the map corresponding to some tangle $T$ with $k$ internal disks, i.e.
$\rho_T(x_1 \otimes x_2 \otimes ... \otimes x_k)$ is equal to $Z(T)$ when the inputs
are $(x_1, x_2, ..., x_k)$.
Then as in \cite{JonB} we have $\rho_T(x_1^* \otimes x_2^* \otimes  ... \otimes x_k^*) = 
\rho_{T'}(x_1 \otimes x_2 \otimes ... \otimes x_k)^*$, where the tangle $T'$ is the mirror image of $T$.
From \cite{JonP}, this means that the BGPA $P_{\Gamma}$ is a planar $*$-algebra, with involution as above.

For ease of later computation, it will be convenient to describe $V_n^{\pm}$ as a set of bounded linear operators on a certain
Hilbert space.

Let $H_n^{\pm}$ be the Hilbert
space with basis $\{x_p\}$ labelled by the paths of length $n$ on $\Gamma$ whose initial vertex
is even ($+$) or odd $(-)$.  These vector spaces are all infinite dimensional when $\Gamma$ is infinite.
We also define the total path Hilbert space $H$ as the direct sum of the $H_n^{\pm}$'s:
$H = \sum_{n, \pm} H_n^{\pm}$.

Each element $A$ of $V_n{\pm}$ naturally defines a linear map on $H_n^{\pm}$.  A loop $L \in l_n^{\pm}$ may
be described as a pair of paths $(\pi, \epsilon)$ of length $n$ whose endpoints are the same.  Then
$\pi$ and $\epsilon$ correspond to basis vectors $x_{\pi}$, $x_{\epsilon}$ in $H_n^{\pm}$, and we take
$\langle A(x_{\epsilon}), x_{\pi} \rangle = A(L)$.  

If $L \in l_n^{\pm}$ is equal to the pair of paths
$(\pi, \epsilon)$, then the reversed loop $L'$ is equal to $(\epsilon, \pi)$.  It follows from the
above definition of the involution that 
$\langle A^*(x_{\epsilon}), x_{\pi} \rangle = 
\overline{\langle A(x_{\pi}), x_{\epsilon} \rangle}$.
In other words, the involution acts on $V_n^{\pm}$ as the adjoint operation for $B(H_n^{\pm})$.

Every path $p$ has a starting vertex $s(p)$ and a terminal vertex $t(p)$.  For each
vertex $v$ of $\Gamma$ we may
define a projection $s_v \in B(H_n^{\pm})$ 
which fixes the closed linear span
of $\{x_p|s(p) = v\}$, and likewise $t_v$ which fixes the closed span of
$\{x_p|t(p) = v\}$.  The $s_v$'s and $t_w$'s form an abelian algebra; in fact the projections
$\{s_v t_w |v, w\in {\rm vertices~of~}\Gamma\}$ are a partition of unity.  Since $\pi$ and
$\epsilon$ above always have the same endpoints, we have each $A \in V_n^{\pm}$ commuting with $s_v$ and $t_v$
for all vertices $v$.  

 This means that $A \in V_n^{\pm}$ may be thought of as a (potentially infinite)
 sum of operators $A_{vw}$ each acting on the subspace $s_v t_w (H_n^{\pm})$. 
Each $s_v t_w$ is of finite rank, and this rank is universally bounded by local finiteness
of $\Gamma$, so these subspaces have bounded dimension. 
Since the components of $A$ are bounded by definition, it follows that $A$ is bounded in norm as
a linear operator on $H$.  

We then have $V_n^{\pm} \subset B(H_n)^{\pm}$, and in fact $V_n^{\pm} \subset \{s_v, t_w | v, w \in
{\rm vertices~of~} \Gamma\} '$.
Each $s_v t_w B(H_n^{\pm})$ is a finite dimensional matrix algebra, with matrix units given by partial isometries
from $x_p$ to $x_q$ where $p$ and $q$ are paths from of length $n$ from $v$ to $w$.  All such partial
isometries are contained in $V_n^{\pm}$, so $s_v t_w B(H_n^{\pm}) \subset V_n^{\pm}$ for all $v, w$.
$\{s_v, t_w\}'$ consists of bounded formal sums of elements of $s_v t_w B(H_n^{\pm})$, so
from the definition of $V_n^{\pm}$ it then follows that $V_n^{\pm} = \{s_v, t_w\}'$ as operators
on $H_n^{\pm}$.  This is a von Neumann algebra by the bicommutant theorem.

We will freely refer to scalar elements of the $V_n^{\pm}$'s; these are just the scalars as
operators on the appropriate Hilbert space.

To assist with computations using the above notation,
we now define a \textit{concatenation operation} $c$ on the $H_n^{\pm}$'s:  

\begin{dfn}
Let $p$ and $q$ be two
paths on the graph $\Gamma$.
Then
$c(x_p,x_q)$ is zero if $t(p) \neq s(q)$, and otherwise is $x_r$ where
$r$ is the path obtained by first following $p$ and then $q$.  This operation extends
linearly and continuously to maps $H_n^{\pm} \otimes H_m^{\pm} \rightarrow H_k^{\pm}$, where
$k = n + m$ and the signs are chosen appropriately. 
\end{dfn}

From the definition of the bases of the $V_n^{\pm}$'s, every concatenation map is surjective.  Furthermore,
these maps are associative: $c(c(x,y),z) = c(x,c(y,z))$. This means we may freely apply $c$ to multiple inputs 
via the inductive definition $c(x_1, x_2, ... , x_n) = c (c (x_1, x_2, ..., x_{n-1}), x_n)$.

We also define an antilinear path reversal operator ${\rm rev}$:

\begin{dfn}
Let $p$ be a path on the graph $\Gamma$.  Let $q$ be the reverse path of $p$, i.e. the same edges and vertices
taken in the opposite order.  Then ${\rm rev}(x_p) = x_q$.  This operation extends antilinearly and continuously
to maps from
$H_n^{\pm}$ to $H_n^{\pm}$ or $H_n^{\mp}$, depending on the value of $n$.
\end{dfn}

Both $c$ and ${\rm rev}$ are bounded, so they extend to the total Hilbert space $H$.  ${\rm rev}$ is an antilinear
involution.

Later, we will be interested in demonstrating that certain maps on the $V_n^{\pm}$'s commute with the action of the planar operad.
To do this it suffices to show that such maps commute with particular tangles that generate the planar operad (see e.g. \cite{BDG2}).
We now describe the action of one set of such generating tangles in terms of the above notation.  All of these actions
may be readily verified directly from the operad definition.

\begin{center}

The multiplication tangle:

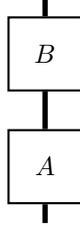
\begin{pspicture}(1.5,3)
\psline[linewidth=2pt](0.75,0)(0.75,0.25)
\psline[linewidth=2pt](0.75,1.25)(0.75,1.75)
\psline[linewidth=2pt](0.75,2.75)(0.75,3)
\psframe(0.25,0.25)(1.25,1.25)
\psframe(0.25,1.75)(1.25,2.75)
\rput(0.75,0.75){$A$}
\rput(0.75,2.25){$B$}

\end{pspicture}

\end{center}

This tangle corresponds to operator multiplication.  The output is $AB$.

\begin{center}

The left embedding tangle $l(A)$:

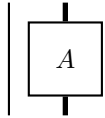
\begin{pspicture}(1.75,1.5)
\psline[linewidth=2pt](1,0)(1,0.25)
\psline[linewidth=2pt](1,1.25)(1,1.5)
\psframe(0.5,0.25)(1.5,1.25)
\rput(1,0.75){$A$}
\psline(0.25,0)(0.25,1.5)
\end{pspicture}

\end{center}

For $A \in V_n^{\pm}$, we have $l(A) \in V_{n+1}^{\mp}$ defined by
$l(A)( c(v,x)) = c(v, A(x))$, where $x \in H_n^{\pm}$ and $v \in H_1^{\mp}$.
This operation is bounded in norm by $e_{max}$.

\begin{center}

The right embedding tangle $r(A)$ :

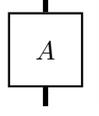
\begin{pspicture}(1.75,1.5)
\psline[linewidth=2pt](0.75,0)(0.75,0.25)
\psline[linewidth=2pt](0.75,1.25)(0.75,1.5)
\psframe(0.25,0.25)(1.25,1.25)
\rput(0.75,0.75){$A$}
\psline(1.5,0)(1.5,1.5)
\end{pspicture}

\end{center}

For $A \in V_n^{\pm}$, we have $r(A) \in V_{n+1}^{\pm}$ defined by
$r(A)(c(x, v)) = c(A(x), v)$, where $x \in H_n^{\pm}$ and $v \in H_1^+$ or
$H_1^-$ depending on the value of $n$.  This operation is norm bounded by $e_{max}$ as well.

Both the left and right embedding operators are strongly continuous; this follows directly from the
definition.

Temperley-Lieb generators ($TL^+$ and $TL^-$): 
\begin{pspicture}(0.9,0.75)
\psarc[fillcolor=lightgray, fillstyle=solid](0.45,0){0.3}{0}{180}
\psarc[fillcolor=lightgray,fillstyle=solid](0.45,0.75){0.3}{180}{360}
\end{pspicture}
and
\begin{pspicture}(0.9,0.75)
\psframe[linestyle=none, linewidth=0pt, fillcolor=lightgray, fillstyle=solid](0,0,)(0.9,0.75)
\psarc[fillcolor=white, fillstyle=solid](0.45,0){0.3}{0}{180}
\psarc[fillcolor=white, fillstyle=solid](0.45,0.75){0.3}{180}{360}
\end{pspicture}

Using the above operad definition, the $TL^+$ element acts on the path basis as follows:

Let $a$ and $b$ be paths of length 1 on $\Gamma$
which start at positive vertices, $c$ and $d$ paths of length 1 which start
at negative vertices, with $x_a$, $x_b$, $x_c$, $x_d$ the corresponding basis vectors in $V_1^+$ and $V_1^-$.
Then 
we should have $\langle TL (c(x_a, x_c)), c(x_b,x_d) \rangle = 0$ unless $a$ is the reverse of $c$ and
$b$ is the reverse of $d$, with additionally $s(a) = s(b)$.  
If these conditions hold, then the inner product should be
$\mu(t(a))\mu(t(c))/\mu(s(a))^2$.

In other words, for each vertex $v \in P_0^+$, let 
$$y_v = \sum_{e | s(e) = v} \mu(t(e)) c(x_e,{\rm rev}( x_e))  $$
where the sum is taken over all paths of length 1 on $\Gamma$.

Then

$$TL(y_v) = \sum_{e_1,e_2 | s(e_{1,2}) = v} 
\mu(t_{e_1}) \mu(t_{e_2})^2 / \mu(v)^2
c(x_{e_2},{\rm rev}(x_{e_2}))$$
$$ = 
\left(
\frac{\sum_{e | s(e)=v} \mu(t(e))^2}{\mu(v)^2}\right)
\left(
\sum_{e | s(e)=v} \mu(t(e)) 
c(x_e,{\rm rev}(x_e))\right)$$
$$ = \frac{\sum_{e | s(e)=v} \mu(t(e))^2}{\mu(v)^2} y_v$$

For each positive vertex $v$,
$y_v$  is an eigenvector for $TL$ with eigenvalue 
$\delta_v = \sum_{e|s(e) = v} \mu(t(e))^2 / \mu(v)^2$, and $TL$ is zero off the
closed linear span of these $y_v$'s.  
We also have a $TL^-$ element.  This is the same diagram as above but with reversed
shading, and is defined by reversing all signs in the above definition.
$y_v$ is defined as above, but now for $v$ being a negative vertex.   Then $TL^-$ has each such
$y_v$ as an eigenvector with eigenvalue $\sum_{e|s(e) = v} \mu(t(e))^2 / \mu(v)^2$, and
 is zero off the closed linear span of these $y_v$'s.

We can now compute left and right capping operators, respectively 
$LC(A): V_n^{\pm} \rightarrow V_{n-1}^{\mp}$ and
 and $RC(A) : V_n^{\pm} \rightarrow V_{n-1}^{\pm}$.

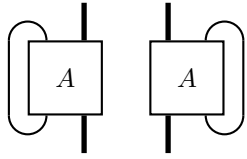
 \begin{center}
 
\begin{pspicture}(1.75,2)
\psframe(0.5,0.5)(1.5,1.5)
\rput(1,1){$A$}
\psline[linewidth=2pt](1.25,0)(1.25,0.5)
\psline[linewidth=2pt](1.25,2)(1.25,1.5)
\psline(0.25,0.5)(0.25,1.5)
\psarc(0.5,0.5){0.25}{180}{360}
\psarc(0.5,1.5){0.25}{0}{180}
\end{pspicture}
\begin{pspicture}(1.75,2)
\psframe(0.25,0.5)(1.25,1.5)
\rput(0.75,1){$A$}
\psline[linewidth=2pt](0.5,0)(0.5,0.5)
\psline[linewidth=2pt](0.5,2)(0.5,1.5)
\psline(1.5,0.5)(1.5,1.5)
\psarc(1.25,0.5){0.25}{180}{360}
\psarc(1.25,1.5){0.25}{0}{180}
\end{pspicture}

\end{center}

These operators may be described in terms of embedding and the TL generators:
 $$r^{n-1}(TL) l(A) r^{n-1}(TL) = r^{n-2}(TL) l^2(LC(A))$$
 $$l^{n-1}(TL) r(A) l^{n-1}(TL) = l^{n-2}(TL) r^2(RC(A))$$

This uniquely defines the capping operation: from the definitio
of left and right embedding $r^{n-2}(TL) l^2(LC(A))$
and $l^{n-2}(TL) l^2(RC( A))$ are each zero only if $LC(A)$, $RC(A)$ respectively are zero.
Note that a graded bounded linear map on the $V_n^{\pm}$'s which fixes the Temperley-Lieb algebra
and commutes with multiplication and embedding also commutes with capping;
for such a map $\omega$ we have
$$l^{n-2}(TL) r^2(RC(\omega(A))) = $$
$$\omega(l^{n-2}(Tl) r^2(RC(A)))=$$
$$l^{n-2}(TL) r^2(\omega(RC(A)))$$ by the above definition, and the same holds for $LC$.

We now describe the interaction of the involution on $V_n^{\pm}$ (defined above) with these tangles.
Since conjugation by the involution corresponds to tangle reflection (see \cite{JonP}), it follows that
$(AB)^* = B^*A^*$, $TL^{\pm}$ are self-adjoint, and the involution commutes with left and right embedding.
This is consistant with the above description of this involution as the adjoint operation on each $B(H_n^{\pm})$.

One important property of the BGPA in \cite{JonB} was the existence of a positive definite $V_0^{\pm}$-valued
sesquilinear form on $V_n^{\pm}$, namely $\langle x,y \rangle = RC^n(y^* x )$ in the above notation.  We would
like this form to be positive definite here as well.  Let $p, q, r, s$ be path basis elements in $V_n^{\pm}$,
and $A, B$ rank one partial isometries from (respectively) $p$ to $q$ and $r$ to $s$.  Then it follows directly
from the operad definition that $\langle A, B \rangle$ is a positive scalar multiple
of a rank one projection in $V_0^+$ if $p=q$ and $r=s$, and is zero otherwise.  It follows that the form
is positive definite on bounded formal sums of such elements, which constitute all of $V_n^{\pm}$.

In order for a BGPA to be useful in a subfactor context, we should have both shaded and unshaded circles
being equal to some scalar $\delta$.  This condition on a planar algebra is called modulus $\delta$ \cite{JonP}.  By 
capping off the single vertical strand $l(1)$, we see this is true when $TL^2 = \delta TL$ with both shadings.  
This occurs when 
$$\delta_v = \sum_{e|s(e)=v} \mu(t(e))^2/\mu(v)^2 $$ is independent of $v$, or 
$$\sum_{e|s(e)=v} \mu(t(e))^2 = \delta \mu(v)^2$$ for all $v$.  Another way of describing this situation is
that the vector with components
$\mu(v)^2$ is an eigenvector for the connection matrix of the graph $\Gamma$, with eigenvalue $\delta$.
Under these circumstances we will say that the spin vector itself (or the BGPA) has modulus $\delta$.

A modulus $\delta$ spin vector is necessarily locally bounded and locally finite, with 
$M \le \sqrt{\delta}$ and $e_{max} \le \delta^2$.  Note that on a finite graph,
there is only one modulus $\delta$ spin vector (up to normalization).  
Its entries are the square roots of the Perron-Frobenius eigenvector
entries for
the inclusion matrix, and $\delta$ is the Perron-Frobenius eigenvalue.  On an infinite locally finite graph, there may be many
modulus $\delta$ spin vectors, with different values of $\delta$.

\section{Automorphisms of BGPAs}

An automorphism of a planar algebra is a graded linear map on the $V_n^{\pm}$'s which commutes
with the entire planar operad.  If the planar algebra has an involution, we will require
the map to commute with the involution as well.  In this section we describe the automorphism
group of an arbitrary BGPA.

These automorphisms of planar algebras may be viewed as a generalization of the automorphisms
of the standard invariant of a subfactor (see \cite{Loi}; c.f. \cite{Sve},\cite{Gup} for specific examples).  
In the BGPA case,
these automorphisms are similar to those of a Jones tower of finite-dimensional $C^*$-algebras,
computed in \cite{Burs}. 

\begin{lem}
Let $\Gamma$ be a bipartite graph, with path Hilbert space $H$ as in section 2.  Let $U$
be a unitary operator on $H$ which respects the grading and commutes with the concatenation operator $c$.  Then
the action of $U$ on $H_0^{\pm}$ is described by a graph automorphism.  Moreover, ${\rm Ad} U$
leaves the BGPA $P_{\Gamma}$ invariant.
\end{lem}
\begin{pf}
Let $x = x_v$ be a standard basis element of $H_0^+$, corresponding to a vertex $v$.
We have $c(x,x) = x$, so $c(U(x), U(x)) = U(x)$ as well by the properties of $U$.
The only elements of $H_0^+$ which have this property are of the form $\sum_{w \in S} x_w$,
where $S$ is some subset of the even vertices of $\Gamma$.  Unitarity of $U$ implies that
in fact $U(x_v) = x_w$ for some $w$.  So $U$ acts by permutation on the even vertices of
$\Gamma$, and likewise on the odd vertices by an identical argument.  Let this permutation
be $\sigma$.

As in section 2, let $s_v$ be the projection onto the basis elements corresponding to
paths starting at $v$,
and $t_v$ the projection onto paths terminating at $v$.  The dimension of $t_v s_w(H_1^+)$
is the number of edges between $v$ and $w$, or zero if they are not adjacent.  Let $n(v,w)$ be this
number of edges.

For any path $p$ from $v$ to $w$ (with corresponding basis element $x_p$), 
we have $c(x_v, x_p,x_w) = x_p$, implying that
$c(x_{\sigma(v)}, U(x_p), x_{\sigma(w)}) = U(x_p)$.  In other words, 
$U t_v s_w (H_1^+) \subset  t_{\sigma(v)} s_{\sigma(w)} (H_1^+)$.  Moreover $U^*$ also
commutes with concatenation and respects grading, so applying the above argument to $U^*$
gives us equality of the above subspaces.  This means that $n(\sigma(v),\sigma(w))=n(v,w)$,
and $\sigma$ is a graph automorphism.

The above also implies that ${\rm Ad} U$ leaves the algebra generated by the $s$'s and $t$'s invariant, since
$U s_vt_w U^* = s_{\sigma(v)}t_{\sigma(w)}$.  Therefore
${\rm Ad} U$ leaves the commutant of this algebra invariant as well.  This commutant is precisely
$V_n^{\pm}$, implying that ${\rm Ad} U$ acts on the planar algebra as desired.

\end{pf} 

We recall that the path reversal  operator ${\rm rev}$ defined in section 2 is an involution
sending $H_1^+$ to $H_1^-$.
From the definition, this operator obeys
$$\langle {\rm rev}(x), {\rm rev}(y) \rangle = \langle y,x \rangle$$
where $\langle \cdot,\cdot \rangle$ is the inner product on $H_1^+$ or $H_1^-$.

\begin{thm}
Let $\Gamma$ be a locally finite bipartite graph with locally bounded spin vector.
Let $H$ be the path Hilbert space on $\Gamma$ as above.
Let $U$ be a unitary element of $B(H)$ which agrees with the grading of $H$
and commutes with the concatenation operator.  Assume further that the restriction
of $U$ to $H_1^{\pm}$ commutes with the path
reversal operator, and that the vertex permutation induced by $U$ 
scales the measure on $H_0$.

Then ${\rm Ad} U$ commutes with the entire planar operad (and adjoint) on $V_n^{\pm}$.

\end{thm}
\begin{pf}
We know from lemma 3.1 that ${\rm Ad} U$ acts on the planar algebra.

Recall $l$ is the left embedding operator.  Let $x$ be in $V_n^{\pm}$, 
$v \in H_1^{\pm}$ and $w \in H_n^{\mp}$.  Then since $U$ commutes with concatentation,
$$U^*l(x) U c(v,w) = U^* l(x) c(U(v),U(w))$$
$$ = U^* c(U(v),xU(w))
= c(v, U^*xU(w))$$
while
$$l(U^* x U)c(v,w) = c(v,U^*xU(w))$$ which is the same.  So
${\rm Ad} U$ commutes with left embedding, and with right embedding as well by a similar argument.

Since $U$ is unitary, ${\rm Ad} U$ commutes with operator multiplication and taking adjoint.  It
remains only to show that ${\rm Ad} U$ commutes with the Temperley-Lieb diagrams, ie with the
elements of the operad representing tangles with no internal disks.  This algebra is generated
by the generators $TL^+$ and $TL^-$
along with multiplication and embedding, so we need to
show that ${\rm Ad} U$ fixes these generators, i.e. that $U$ commutes with them. 

Let $v$ and $w$ be two adjacent vertices of $\Gamma$.  Let $H_1^{vw}$ be the subspace
of $H_1$ spannned by paths from $v$ to $w$, and $H_2^{vw}$ the subspace of $H_2$
spanned by paths from $v$ to $w$ and back to $v$.  Let $\{a_i\}$ be any orthonormal basis
for $H_1^{vw}$.  Take $y_{vw} = \sum_i c(a_i , {\rm rev}(a_i))$, using the reversal operator
defined above.  Because of the interaction of ${\rm rev}$ with inner product described
above, the inner product of $y_{vw}$ with any element of the form $c(v, {\rm rev}(v))$
is equal to the squared norm of $v$ regardless of which basis is chosen.  Such elements span
$H_2^{vw}$, so $y_{vw}$ is independent of the choice of specific basis.  This is true
if $v$ is odd or even.

Let $\sigma$ be the permutation action of $U$ on vertices of $\Gamma$.
As in lemma 3.1, we have $U(H_1^{vw}) = H_1^{\sigma(v) \sigma(w)}$, and so
$U$ maps an orthonormal basis for the first vector space to one for the second. 
This means that
$U( y_{vw} ) = y_{\sigma(v)\sigma(w)}$.

We recall from section 2 that the $TL$ generators in $B(H_2^{\pm})$ each leave the closed linear span
of certain vectors $\{x_v \}$ invariant, where $v$ is taken from the set of positive or negative
vertices depending on the sign
of the generator, and are zero off the span of these vectors.  We defined
$$y_v = \sum_{e|s(e)=v} \mu(t(e)) c(x_e, {\rm rev}(x_e))$$ for each vertex $v$
but we 
can also write $$y_v = \sum_{w|n(v,w) \neq 0} \mu(w) y_{vw}$$ using the notation above.

From lemma 3.1 $\sigma$ is a graph automorphism of $\Gamma$.  By assumption,
$\sigma$ scales the spin vector by a fixed constant $\lambda$.  Therefore 
$$U(\mu(w) y_{vw}) =
\lambda^{-1}(\mu(\sigma(w)) y_{\sigma(v)\sigma(w)}$$
Since $\sigma$ is a graph automorphism, it maps the set of vertices adjacent
to $v$ to the set of vertices adjacent to 
$\sigma(v)$, and we have as well
$$U \left( \sum_{w|n(w,v)\neq 0} \mu(w) y_{vw} \right) =
\sum_{w|n(w,\sigma(v)) \neq 0}\lambda^{-1} \mu(w) y_{\sigma(v)w}$$

So $U$ maps one standard basis vector of the subspace acted on by either $TL$ generator to $\lambda^{-1}$
times another such vector.  

We have $TL(y_v) = \delta_v y_v$, where 
$\delta_v=\sum_{e|s(e) = v} \mu(t(e))^2/\mu(v)^2$.  From the properties of $\sigma$, 
$$\delta_{\sigma(v)} = 
\sum_{e|s(e) = \sigma(v)} \mu(t(e))^2/\mu(\sigma(v))^2
=  
\sum_{e|s(e) =v} \mu(\sigma(t(e)))^2/\mu(\sigma(v))^2$$
$$
=\sum_{e|s(e) =v} \lambda^2 \mu(t(e))^2/\lambda^2\mu(\sigma(v))^2
=\sum_{e|s(e) =v} \mu(t(e))^2/\mu(v)^2 = \delta_v
$$
So $U$ leaves each eigenspace of both $TL$ generators invariant.

This means
that $U$ commutes with the $TL$ generators.  So $U$ commutes with a set of generating
tangles for the planar operad, and hence with the entire operad.

\end{pf}

Now we describe two classes of linear maps on a BGPA which meet all the above conditions.

\begin{dfn}
Let $\Gamma$ be a bipartite graph with spin vector $\mu$.  Let $\kappa$ be a permutation
of the vertices of $\Gamma$, preserving connection numbers and parity and 
preserving or scaling spin.  
Label each $n$-fold multiple edge by $\{1,...,n\}$.  We may then extend $\kappa$ to the edges of
$\Gamma$
by asserting that it preserves this numbering. 
Then $\kappa$ gives rise to a permutation of the paths on $\Gamma$, 
and therefore a map $U$ on the path Hilbert space $H$.  Then ${\rm Ad} U$ is the 
\textit{graph automorphism
operator} associated with $\kappa$.
\end{dfn}

\begin{lem}
Graph automorphism operators are automorphisms of the BGPA.
\end{lem}
\begin{pf}
Let $Ad U$ be a graph automorphism operator as above.
It follows directly from the definition that $Ad U$ commutes with path reversal and concatenation,
and agrees with the grading.  We have also assumed that the underlying graph automorphism is 
trace scaling.  So all of the conditions of theorem 3.1 are satisfied, and $Ad U$ is a planar
algebra automorphism.
\end{pf}

\begin{dfn}

Let $O$ be an element of $V_1^+$ which is unitary as
an operator acting on $H_1^+$.  Let $O' = {\rm rev} \circ O \circ {\rm rev}$, acting
on $H_1^-$.  For $1 \le i \le n$, let $p_i$ be a path basis element of $H_1^+$ (for $i$ even)
or $H_1^-$ (for $i$ odd).  Then let
$$U(c(p_1,p_2,...,p_n)) = c(O(p_1), O'(p_2), O(p_3) ...)$$
and $U$ extends to a unique bounded linear operator on $H_n^+$.  Define
$U$ similarly on $H_n^-$ as the extension of the map
$$U(c(p_1,p_2,...,p_n)) = c(O'(p_1), O(p_2), O'(p_3) ...)$$
where the $p_i$'s are again in $H_1^+$ or $H_1^-$ as appropriate.

Then $Ad U$ is the \textit{multiplication operator} associated with $O$.
\end{dfn}

If $O$ acts nontrivially on only one of the subspaces $s_v t_w (H_1^+)$, while
leaving all others fixed, then we will call it a \textit{basic} multiplication operator associated
to the
vertex pair $\{v,w\}$.  Every multiplication operator is a product of basic multiplication operators,
and basic multiplication operators associated to different vertex pairs commute with each other.
A multiplication operator is \textit{scalar} if the restriction to $B(s_v t_w (H_1^+))$ for each 
$v$, $w$
is a scalar multiple of the identity.  The scalar multiplication operators are the center
of the multiplication operator group.  If the graph has no multiple edges then every
multiplication operator is scalar.

\begin{lem}
Multiplication operators are automorphisms of the BGPA.
\end{lem}
\begin{pf}
Let ${\rm Ad} U$ be a multiplication operator as above.  From the definition, it commutes with
concatenation and path reversal and respects grading.  Since $O$ comes from an element
of $V_1^+$, it commutes with $s_v$ and $t_v$ for vertices $v$. Therefore the associated graph
automorphism $\sigma$ is trivial, and preserves the trace.  So all the conditions of theorem
3.1 are satisfied, and $Ad U$ is a planar algebra isomorphism.
\end{pf}

Direct computation shows that conjugating a basic multiplication operator by a graph automorphism
operator
produces a basic multiplication operator associated to different vertex pair.  So the group of automorphisms
generated by these two types of operators has a crossed product structure: the subgroup
of multiplication operators is normal.

Now we will show that the two types of operators described above in fact generate the entire
automorphism group of a BGPA.

\begin{lem}
Any automorphism of a BGPA is strongly continuous.
\end{lem}
\begin{pf}
Let $\alpha$ be an automorphism of a BGPA.

From the definition in section 3.2,
we have $V_n^{\pm} \subset B(H_n^{\pm})$, where $V_n^{\pm}$ is a (potentially infinite)
direct sum of type $I$ factors.  
Every automorphism of such a von Neumann algebra can be written ${\rm Ad} U$, where $U$ is a unitary
element of  $B(H_n^{\pm})$. Since $\alpha|_{V_n^{\pm}}$ commutes with multiplication and 
involution, this restriction of $\alpha$ is a von Neumann algebra isomorphism, and may be written
${\rm Ad} U_n^{\pm}$ for unitary $U \in B(H_n^{\pm})$.  Summing all the $U_n^{\pm}$s provides a unitary
operator on the total Hilbert space $H$
whose adjoint action on the graded vector space $V_n^{\pm}$ agrees with $\alpha$. 

Multiplication is strongly continuous.

\end{pf}

This means that two automorphisms of a BGPA are equal if and only if they agree on loops,
since the loops span a strongly dense set in each $V_n^{\pm}$.

\begin{lem}
Let $\alpha$ be an automorphism of a BGPA.  Then there is a graph automorphism operator
$\beta$ such that $\alpha$ and $\beta$ agree on $V_0^{\pm}$.
\end{lem}
\begin{pf}

Let $p_v$, $q_w$ be the atomic projections in $V_0^+$, $V_0^-$ associated with even and odd vertices
$v$ and $w$.  Since $\alpha$ is a BGPA automorphism, it must send $p_v$ to some other atomic projection
$p_{\sigma(v)}$ 
in $V_0^+$, and likewise for $q_w$.  Therefore $\alpha$ induces a permutation $\sigma$
on the vertices of $\Gamma$.

To see that this permutation is a graph automorphism, note that $l(p_v) r(q_w)$ is a minimal central
projection of $V_1^+$ (or zero), and the dimension of $l(p_v) r(q_w) V_1^+$ is the square of the number of edges
between $v$ and $w$.  This dimension is preserved by $\alpha$, i.e. it agrees with the dimension of
$l(p_{\sigma(v)}) r(q_{\sigma(w)}) V_1^+$.  Therefore $(v,w)$ and $(\sigma(v),\sigma(w))$ have the same
number of edges between them, i.e. $n(v,w) = n(\sigma(v),\sigma(w))$.

Next we must show that $\sigma$ preserves or scales the trace.  For this we note that $p_v$ with a circle
around it evaluates to 
$$\sum_{w|n(v,w) \neq 0} q_w \frac{\mu(v)}{\mu(w)}$$
  Since $\alpha$
commutes with the tangle, we must have 
$$\mu(v)/\mu(w) = \mu(\sigma(v))\mu(\sigma(w))$$ for all adjacent
$v$, $w$, implying the desired result.

From the above description of graph automorphism operators, there is a such an operator
$\beta$ whose induced permutation action on the vertices of $\Gamma$ is the same as $\alpha$'s.  These
two operators then agree on $V_0^{\pm}$.

\end{pf}

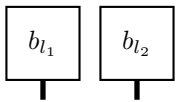
\begin{lem}
Let $\alpha$ be an automorphism of a BGPA which acts trivially on $V_0^\pm$.  Then there
is a multiplication operator $\beta$ such that $\beta$ and $\alpha$ agree on $V_1^{\pm}$.
\end{lem}
\begin{pf}
Since $\alpha$ acts trivially on $V_0^{\pm}$, it fixes all elements of the form $r(p_v) l(q_w)$ in 
$V_1^+$.  These elements are the center of $V_1^+$ taken as a von Neumann algebra, 
so $\alpha$ acts as an inner automorphism on $V_1^+$.  There
is a multiplication operator $\beta$ whose action on $V_1^+$ is any desired inner automorphism.  Then
$\alpha$ and $\beta$ agree on $V_1^+$.  Both of these automorphisms commute with the half rotation

\begin{center}

\begin{pspicture}(2.5,2.5)

\psframe(0.75,0.75)(1.75,1.75)
\rput(1.25,1.25){$A$}
\psarc(0.75,0.75){0.5}{180}{360}
\psarc(1.75,1.75){0.5}{0}{180}
\psline(0.25,0.75)(0.25,2.5)
\psline(2.25,1.75)(2.25,0)
\end{pspicture}

\end{center}

which is a bijective map from $V_1^+$ to $V_1^-$, so they agree on $V_1^-$ as well.
\end{pf}

\begin{lem}
Let $\alpha$ be an automorphism of a BGPA which acts trivially on $V_1^{\pm}$.  Then $\alpha$
is a scalar multiplication operator.
\end{lem}
\begin{pf}

From the proof of lemma 3.4, we can write the action of $\alpha$ on each $V_n^{\pm}$ as
${\rm Ad} U_n^{\pm}$, where $U_n^{\pm}$ is a unitary in $B(H_n^{\pm})$.

First note that for any path $p$ on the graph $\Gamma$, with corresponding basis
vector $x_p \in H_n^{\pm}$, there is a rank one projection onto $x_p$ contained
in $V_n^{\pm}$.  This projection may be written as the product
$$r^n(p_{e_1}) r^{n-1}l( p_{e_2}) ... r l^{n-1} (p_{e_{n-1}})l^n (p_{e_n})$$
where $e_i$ is the $i$th edge of $p$ and $p_{e_i}$ is the rank one projection in $V_1^+$ or $V_1^-$ onto
the vector $x_{e_i}$ corresponding to the path $e_i$.  Since $\alpha$ acts trivially on $V_1^{\pm}$, it fixes this projection, and
$U_n^{\pm}$ must therefore map each $x_p$ to a scalar multiple of itself.

Now let $b_l \in V_n^{\pm}$ be a rank one partial isometry corresponding to a loop $l \in l_n^{\pm}$.
We have
$b_l = x_p b_l x_q$ for certain rank one projections $x_p$,
$x_q$ as above, implying that 
$$\alpha(b_l) = \alpha (x_p b_l x_q) = x_p \alpha(b_l) x_q$$
The only way this can be true is if $\alpha$ sends $b_l$ to a scalar multiple of itself as well.

In other words, every loop is an eigenvector of $\alpha$, and $\alpha$ induces a map $\rho$ from the set of loops
to the complex scalars of modulus 1.  

Since $\alpha$ commutes with the half rotation,
$\rho(l)$ is independent of the basepoint of $l$.  Because $\alpha$ fixes $V_1^{\pm}$, $\rho(l) = 1$ for any loop
$l$ of length 2.  Finally, since $\alpha$ commutes with this diagram 

\begin{center}

\begin{pspicture}(2.75,1.25)
\psline[linewidth=2pt](0.75,0)(0.75,0.25)
\psframe(0.25,0.25)(1.25,1.25)
\rput(0.75,0.75){$b_{l_1}$}

\psline[linewidth=2pt](2,0)(2,0.25)
\psframe(1.5,0.25)(2.5,1.25)
\rput(2,0.75){$b_{l_2}$}

\end{pspicture}

\end{center}

if $l_3$ is the concatenation of the loops $l_1$ and $l_2$ then we have
$\rho (l_3) = \rho(l_1) \rho(l_2)$.  

Putting this together we find that $\rho$ is necessarily
a 1-dimensional representation of the fundamental group of $\Gamma$.  But any such representation may be obtained
from a scalar multiplication operator.  We may always find a set of free generators $\{l_1, l_2, ...\}$ for the fundamental group 
with the property that each generator $l_i$ contains some edge $e_i$ which does not appear in any other loop.  Then 
a basic scalar multiplication operator with value $\lambda$ associated to the endpoints of
some $e_k$ corresponds to the representation of $\pi_1(\Gamma)$ sending
$l_k$ to $\lambda$ (or possible $\overline{\lambda}$, depending on the direction of $l_k$) and all other generators
to 1.  All other representations of $\pi_1$ may be obtained similarly from basic scalar multiplication operators 
associated to various $e_i$'s.

This implies that there is a scalar
multiplication operator which agrees with $\alpha$ on loops.  
Since the loops span a strongly dense subset of each
$V_n^{\pm}$, and automorphisms of a BGPA are strongly continuous by lemma 3.4, it follows that $\alpha$ itself
is a scalar multiplication operator.

\end{pf}

Since scalar multiplication operators themselves act trivially on $V_1^{\pm}$, this lemma in fact shows that the scalar
multiplication operators of a BGPA are isomorphic to the 1-dimensional representations of the fundamental group of the
graph.

\begin{thm}
Let $P_{\Gamma}$ be a BGPA, with multiplication operators $E$ and graph automorphism operators $A$.
Let $\alpha$ be an automorphism of $P_{\Gamma}$.  Then $\alpha = ae$ for some 
$a \in A$, $e \in E$.
\end{thm}
\begin{pf}
This follows from the proceeding lemmas.  There is $\beta_1 \in A$ such
that $\beta_1^{-1} \alpha$ acts trivially on $V_0^{\pm}$.  There is $\beta_2 \in E$ such
that $\beta_2^{-1} \beta_1^{-1} \alpha$ acts trivially on $V_1^{\pm}$.  So
$\beta_2^{-1} \beta_1^{-1} \alpha$ is a scalar multiplication operator $\beta_3$, and
$\alpha = \beta_1 \beta_2 \beta_3$ with $\beta_1 \in A$ and $\beta_2 \beta_3 \in E$.
\end{pf}

Since conjugation by graph automorphisms leaves the multiplication operator group invariant, we may in fact
write ${\rm Aut} P_{\Gamma} = E \rtimes A$ with notation as above.

\section{Planar fixed point subfactors}

A subfactor planar algebra is a planar algebra with the following additional properties ~\cite{JonP}:

\begin{itemize}
\item
${\rm dim} V_0^{\pm} = 1$
\item
${\rm dim} V_n^{\pm} < \infty~\forall n$
\item
Spherical: Since $V_0^+ = V_0^- = \mathbb{C}$, we may equate these vector spaces with the scalars.  A planar
algebra is \textit{spherical} if
for any $A \in V_1^+$, the left and right caps $LC(A)$ and $RC(A)$ agree.
\item
Involution: There is an antilinear isometry on each $V_n^{\pm}$ which interacts with tangles as reflection.
\item
Positive definiteness: Involution gives us a scalar sesquilinear form, namely $\langle x,y \rangle = RC^n(y^* x)$
where the multiplication and right capping tangles are as in section 2.
This form should be
positive definite.
\end{itemize}

The standard invariant of any finite index extremal $II_1$ subfactor may be described as a subfactor
planar algebra (\cite{JonP},\cite{Pop}).  We take $V_n^+ = M_0' \cap M_n$, $V_n^- = M_1' \cap M_{n+1}$, and the operad
definition is given in \cite{JonP}.

Conversely, if $P$ is a subfactor planar algebra, then there exists a finite index extremal
$II_1$ subfactor such that the standard invariant of this subfactor is $P$ (\cite{Pop}, c.f. \cite{JonP},\cite{GJS}, 
\cite{KS},\cite{JSW}).

For finite graphs with the correct spin vector, any sufficiently small planar subalgebra
is of subfactor type.

\begin{lem}
Let $P_{\Gamma}$ be a bipartite graph planar algebra with spin vector $\mu$.  
Let $x$ be in $V_1^+$; let $x_l$ represent the element of $V_0^-$ obtained by capping off to the left,
and $x_r \in V_0^+$ from capping off to the right.  Suppose both $x_l$ and $x_r$ are scalars. 
Then there is some constant $\alpha$, independent of $x$,
such that $x_l = \alpha x_r$.
\end{lem}
\begin{pf}
From ~\cite{JonB}, there is a partition function defined on $V_0^+$ and $V_0-$ as the linear extension of
$x_v \rightarrow \mu^4(v)$, and this function has the same value on $x_l$ and $x_r$.

If $x_l$ and $x_r$ are scalars, this means $x_l = \lambda_1 \sum_v x_v$, $x_r = \lambda_2 \sum_v x_v$.

The partition function on $x_1$ is $\lambda_1 \sum_v \mu^4(x_v)$, and on $x_r$ is $\lambda_2 \sum_v \mu^4 (x_v)$.
Since these are the same we must have $\lambda_1 = \lambda_2 \frac{\sum_v \mu^4(x_v)}{\sum_v \mu^4(x_v)}$.  In other
words $\alpha = \frac{\sum_v \mu^4(x_v)}{\sum_v \mu^4(x_v)}$ and $x_l = \alpha x_r$ for all $x$ such that
$x_l$ and $x_r$ are scalars.

\end{pf}

\begin{thm}
Let $P_{\Gamma}$ be a finite bipartite graph planar algebra whose spin vector $\mu$ has modulus $\delta$.
Let $X$ be a planar $*$-subalgebra of $P_{\Gamma}$ such that $X \cap V_0^+ = X \cap V_0^- = \mathbb{C}$.
Then $X$ is a subfactor planar algebra.
\end{thm}
\begin{pf}
The BGPA has an involution, giving rise to a positive definite $V_0^+$-valued sesquilinear form
(see section 2).  The involution and form are inherited by $X$, and the restriction of
the form to $X$ is scalar valued since $X \cap V_0^{\pm} = \mathbb{C}$.  Since $\Gamma$ is finite,
each $V_n^{\pm}$ is finite dimensional, and so this is true of their intersections with $X$ as well.

To show that $X$ is of subfactor type, it remains only to demonstrate sphericality.  Let $x \in X$
be an element of $V_1^+$.  Capping off to left or right produces scalars, since $X \cap V_0^{\pm}$ is scalar.  
Therefore the conditions of the lemma above are satisfied.  Because $\mu$ has modulus $\delta$,
shaded and unshaded circles represent the same scalar.  These diagrams
are the left and right caps of a single vertical strand, so the constant $\alpha$ in the lemma
is equal to 1.  It follows that $x_l$ and $x_r$ are equal as scalars, and $X$ is spherical.
\end{pf}

A small subalgebra of an infinite BGPA still corresponds to a subfactor, but the subfactor need not be extremal.

\begin{thm}
Let $P_{\Gamma}$ be a locally finite bipartite graph planar algebra.  Let $X \subset P_{\Gamma}$ be a
spherical planar $*$-subalgebra with $X \cap V_0^{\pm} = \mathbb{C}$.  Then $X \cap V_n^{\pm}$ is finite dimensional.
\end {thm}
\begin{pf}

Let $p_v \in V_0^+$ be a minimal projection corresponding to some even vertex $v$, and $q = r^n(p_v)$.
Take $x \in X \cap V_n^+$, and
let $RC^n$ be the diagram consisting of capping off all strands to the right (this is the form from
section 2).  Then
$RC^n(x)$ is a scalar from the properties of $X$, and $RC^n(q x) = qRC^n(x)$ by isotopy invariance.
$q$ commutes with $x$ and $x^*$ with respect to the usual multiplication tangle.
If $q x = 0$, then $qx^* x = 0$, giving
$RC^n(q x^* x) = 0 = q RC^n(x^* x)$.  Since $RC^n(x^*x)$ is a scalar, this means that $px = 0$ implies
$RC^n(x^* x) = 0$.  But $RC^n$ gives a positive definite form on $V_n^+$, so $x$ itself is zero in this case.

This means that the map $x \rightarrow qx$ is injective on $X \cap V_n^+$.  But
$q V_n^+$ has basis labelled by loops of length $2n$ which start and end  at $v$.  By local
finiteness of $\Gamma$, this set is finite.  Therefore $p V_n^+$ is finite dimensional, and
$X \cap V_n^+$ is as well.

The same argument shows that $X \cap V_n^-$ is finite dimensional.
\end{pf}

Burns described a class of \textit{rigid planar $C^*$-algebras} in \cite{Burns}, generalizing Jones' definition of a SPA.  
Every rigid planar $C^*$-algebra is the standard invariant of a finite index $II_1$ subfactor, but this subfactor
need not be extremal.  From Burns' definition, a planar algebra having all the characteristics of a SPA except sphericality is
a rigid planar $C^*$-algebra.

Let $P_{\Gamma}$ be a modulus $\delta$ BGPA coming from an infinite finite graph, 
and $X \subset P_{\Gamma}$ a planar subalgebra with ${\rm dim } (X \cap V_0^{\pm})=1$.  
The above lemma tells us that $X \cap V_n^{\pm}$ is finite dimensional.  The BGPA has a positive
definite sesquilinear form and an involution with the right properties, which are inherited by $X$.  This means that $X$ 
is a rigid planar $C^*$-algebra.  If $X$ is additionally spherical, then it is an SPA.  This is automatic when $X$ is 
irreducible, but is not true in general.

This means that there is a finite index $II_1$ subfactor whose standard invariant is $X$.  The subfactor will be extremal
if and only if $X$ is spherical.

\section{Examples}

\subsection{Introduction}

We can construct a wide range of subfactor planar algebras under this fixed point technique.  In order
for the fixed points $P_{\Gamma}^G$ to be an SPA, we need it to have 1-dimensional intersection
with $V_0^+$ and $V_0^-$, which is equivalent to
having the graph automorphism part  of $G$ act transitively on both positive and negative vertices.  We should
also check sphericality for non-irreducible infinite graph examples.
Some SPAS thus are the standard invariants of previously known subfactors, while others seem to be 
previously unclassified.  A few such examples are described below.

\begin{dfn}
If $P_{\Gamma}$ is a BGPA, and $G \subset {\rm Aut} P_{\Gamma}$ with $P_{\Gamma}^G$ an SPA,
then the corresponding subfactor is a planar fixed point subfactor.
\end{dfn}

\subsection{Group-subgroup subfactors}
A specific example of the planar fixed point construction is found in \cite{Gup}.  In this
paper, Gupta starts 
with the BGPA on the graph with $n$ odd vertices all connected to one even vertex.
Then $G$ is some group acting by permutation on the set of odd vertices, and $H$ is
the subgroup which fixes some specified vertex.  Gupta shows that the fixed points of
the BGPA by this action constitute a subfactor planar algebra, and that this is in fact
the standard invariant of the group-subgroup subfactor (see e.g. \cite{KY})
corresponding to the inclusion 
$H \subset G$, namely $M^G \subset M^H$ for some outer action of some finite group $G$ on 
a $II_1$ factor $M$.

\subsection{Wassermann subfactors}

Here we let $\Gamma$ be the graph with two vertices connected by an $n$-fold multiple
edge.  Let $G$ be any compact subgroup of  the unitaries $M_n(\mathbb{C})$.  Then $G$
may be embedded in the multiplication operators on this graph.  The fixed points of this $G$-action
 are of subfactor type.  They are identical to the standard invariant of the Wasserman
subfactor 
$$(\overline{1\otimes M_n \otimes M_n \otimes ...}^{st})^G \subset (\overline{M_n \otimes M_n \otimes M_n \otimes ...}^{st})^G$$
where $G$ acts pointwise on the tensor products (see \cite{Was}). 

\subsection{Diagonal subfactors}
Let $G$ be a finitely generated group of outer automorphisms of a $II_1$ factor $M$.
Let $G$ have distinct generators $\{g_1, ... ,g_n\}$, and take $\Gamma$ to be the graph which has
one odd and one even vertex for each element of $G$.  Two vertices $v_x^+$ and $v_y^-$
are connected by a single edge if $y =  x h$, for $h \in \{1,g_1,...,g_n\}$.   The spin 
vector of this graph is 1 at every vertex; it has modulus $n+1$.

Then $G$
acts on the graph by left translation: $\alpha_x(v_g^{\pm}) = v_{xg}$.  Associating each
such graph automorphism with the corresponding graph automorphism operator on $P_{\Gamma}$
gives an action of $G$ on $P_{\Gamma}$.  

This action is transitive on even and odd vertices.  When the graph is infinite, sphericality
may be directly verified: $V_1^+$  has dimension $n+1$, and each minimal
projection has left and right trace $1/(n+1)$.  Therefore $P_{\Gamma}^G$ is of subfactor type.

Since vertices are labelled by group elements, loops of length $2n$
in this graph may be written as a list of vertices
$$x-xa_1-xa_1a_2^{-1}xa_1a_2^{-1}a_3...-xa_1a_2^{-1}...a_{2n}^{-1} = x$$
where each $a_i$ comes from the set $\{1, g_1, ..., g_n\}$.  So we may think of this loop
as a starting point $x$ along with a list of generators $(a_1, a_2, ...)$
whose alternating product is the identity in $G$.

The group action moves the base point (via left translation) while keeping the generator
list invariant.  It follows that these generator lists label a basis for the intersection of $P_{\Gamma}^G$
with $V_n^{\pm}$.

This basis is precisely that described in \cite{BDG1} for the planar algebra of the
diagonal subfactor (see \cite{PoF},\cite{BiE}).  It may be verified that the operations of left and right embedding,
involution, and multiplication on $P_{\Gamma}^G$ agree with the planar algebra of \cite{BDG1}, and
the Temperley-Lieb algebra embeds in the same way,
so these planar algebras are isomorphic.  It follows that the planar algebra constructed
in this way is that of a diagonal subfactor without cocycle.

\subsection{Bisch-Haagerup subfactors}
Let $G$ be a group of outer automorphisms of a $II_1$ factor $M$, generated by finite subgroups
$H$ and $K$.  For simplicity we require $H \cap K = \{1\}$.
 
Let $\Gamma$ be the graph which has one even vertex for each $H$ right coset in $G$,
and one odd vertex for each $K$ right coset.  Then the edges of $\Gamma$ are labelled by
group elements.  The endpoints of each edge $e_g$ are the 
even vertex $v_{gH}$ and the odd vertex $v_{gK}$.  The spin vector has value $H$ on each
even vertex and value $K$ on each odd vertex, and it has modulus $\sqrt{|H||K|}$.
 
We may write a loop of length $2n$ as a list of edges:
$$e_{x_1} - e_{x_2} ... - e_{x_{2n}} - e_{1}$$
Assume first that the loop begins at a positive vertex.
For the path to be connected, $e_{x_1}$ and $e_{x_2}$ must share a vertex, so
$x_1$ and $x_2$ are in the same $K$-coset and $x_2 = x_1 k_1$.  Likewise $x_2$ and $x_3$
are in the same $H$-coset.  So we can view this loop as a starting point $x$ along with
a list of alternating elements of $k$ and $h$, with the restriction that 
$k_1h_1k_2h_2...k_nh_n$ is equal to the identity in $G$.

  Again $G$ acts on the graph by left translation,
and this gives a $G$-action on $P_{\Gamma}$.  This action shifts the basepoint of loops
while leaving the list of $k_i$'s and $h_i$'s the same.  So we may identify a basis for
$P_{\Gamma}^G \cap V_n^+$, namely $n$-tuples of elements of $K$ and $H$ obeying
$k_1h_1k_2h_2...k_nh_n=1$; the basis for $V_n^-$ is obtained by reversing the roles
of $H$ and $K$.  This labelling set is the same as that of the planar algebra
of \cite{BDG2} (c.f. \cite{BH}); again it may be shown that the two planar algebras are isomorphic.  It follows
that this fixed point planar algebra is the standard invariant of the Bisch-Haagerup subfactor
$M^H \subset M \rtimes K$.
 
\subsection{The cube graph}

Let $\Gamma$ be the graph of a cube, and $P_{\Gamma}$ the corresponding BGPA. 
We obtain several subfactor planar algebras by taking fixed points under various group
actions.  

We note that the automorphism group of the bipartite graph is $S_4$, since each
such automorphism may be described uniquely as a certain permutation of the even vertices.
The only modulus $\delta$ spin vector assigns weight 3 to every vertex; 
$\delta = 3$, so every subfactor planar subalgebra produces an index 9
subfactor.  We can write down a biprojection which is invariant under every automorphism,
so there is always an index 3 intermediate subfactor.

We mention some of the possibilities described above. 
If $G$ is the subgroup generated by $(12)(34)$ and $(13)(24)$, then $P_\Gamma^G$ is the
planar algebra of the diagonal subfactor; the group generators are each order 2, and the
group is $Z_2^2$.  If $G$ is $A_4$, then $P_\Gamma^G$ is the Bisch-Haagerup subfactor 
$M^{Z_3} \subset M \rtimes Z_3$, where the two order-3 automorphisms of $M$ together
generate the group $A_4$.

If we take $G = S_4$, then $P_{\Gamma}^G$ is some other subfactor planar algebra.  Its principal
graph may be directly computed from the group action:

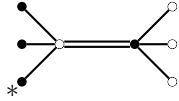
\begin{center}
\begin{pspicture}(2.5,1.5)
\psset{dotstyle=*}
\dotnode(0.25,0.25){A1}
\dotnode(0.25,0.75){A2}
\dotnode(0.25,1.25){A3}
\dotnode(1.75,0.75){C}

\psset{dotstyle=o}
\dotnode(0.75,0.75){B}
\dotnode(2.25,0.25){D1}
\dotnode(2.25,0.75){D2}
\dotnode(2.25,1.25){D3}

\ncline{A1}{B}
\ncline{A2}{B}
\ncline{A3}{B}
\ncline{D1}{C}
\ncline{D2}{C}
\ncline{D3}{C}

\rput(0.125,0.125){$*$}
\psset{offset=1pt}
\ncline{B}{C}
\ncline{C}{B}
\end{pspicture}
\end{center}

The dual principal graph is the same.
It does not appear to be of any previously categorized type, although we have tentatively identified
it as a composition of two group-subgroup subfactors.

Finally, we may take $G = {\rm Aut} P_{\Gamma}$.  Since the multiplication operator group of $P_{\Gamma}$ is
$(S^1)^5$, this group is infinite, and $P_{\Gamma}^G$ is infinite depth.

\subsection{The degree (3,2) tree graph}

Let $\Gamma$ be the graph which branches twice at each even vertex and three times at each odd
vertex.  The operator norm of this graph is $4 \sqrt{2}$.   There are no multiplication operators on
$P_{\Gamma}$, so all automorphisms will come from graph automorphisms.

First assign weight $2$ to every even vertex and $3$ to every odd vertex.  This spin vector
has modulus $\sqrt{6}$.  Every automorphism,
of $\Gamma$ here gives rise to an automorphism of $P_{\Gamma}$.  Taking $G = Aut(P_{\Gamma})$,
we obtain a subfactor planar algebra $P_{\Gamma}^G$.  This subfactor is irreducible, hence
automatically spherical, and is non-amenable; we conjecture
it is obtained as a composition of two group-subgroup subfactors, where the groups involved are
variations of the Grigorchuk lamplighter group.

We may also obtain a transitive subgroup as follows: color the edges of the graph with three colors so that
no vertex contacts two edges of the same color, and then consider all automorphisms
which leave the coloring invariant or permute the colors.  This group is $Z_3 * Z_2$, and
the resulting subfactor is the index 6 Bisch-Haagerup subfactor described in \cite{BH} corresponding to this group.

Now we describe a non-irreducible example.  We color the edges of the graph red and blue so that
each even vertex contacts a red and blue edge, and each odd vertex contacts two reds and
a blue.  We now consider $G$ to be the group of color preserving automorphisms.  It may
be seen that this group is transitive on odd and even vertices, but the fixed points have 2-dimensional
intersection with $V_1^+$ space, corresponding to the 2 $G$-orbits of edges.  With the above spin vector, sphericality
fails: the 'blue edge' element has left and right traces of 1/3 and 1/2.

To find a subfactor planar subalgebra of this BGPA, we will need a different spin vector.  We take
$\mu(v_0) = 1$ for some arbitary base vertex $v_0$, and then define $\mu$ elsewhere
so that for adjacent vertices $x$ even and $y$ odd, we have $\mu(y)=\mu(x)$ if the $x-y$ edge
is red, but $\mu(y) = 2^{1/4}\mu(x)$ if the edge is blue.  This trace has modulus $1 + \sqrt{2}$.
It may be seen that any color preserving automorphism of $\Gamma$ will multiply $\mu$ by a
constant factor, hence providing an automorphism of $P_\Gamma$.  With this choice of spin vector and
$G$ the color preserving graph automorphisms,
the left and right traces of the 'blue' element are both $\sqrt{2}-1$, and sphericality holds.
So from section 5.4, $P_{\Gamma}^G$ is a subfactor planar algebra.

The resulting subfactor is infinite depth non-amenable.  
Its index is 
$(1 + \sqrt{2})^2 = 
3 + 2\sqrt{2}$, 
and it does not have any intermediate subfactors.

From \cite{PP}, this is the minimum index for an extremal non-irreducible subfactor.  The above construction
is a new way of getting such a subfactor.  Any group which is transitive on even and odd vertices
but the same two orbits of edges as above will also produce a subfactor with this list of properties,
so we actually have many such examples.  For example, we might partition the red edges into 'red'
and 'white' so that each odd vertex contacts one vertex of each color; $G'$ might then be the group
of graph
automorphisms which either preserve color or swap the colors red and white.

There are also spin vectors of any modulus greater than $3 + 2\sqrt{2}$ on this graph
which are preserved or scaled by every element of $G$.  Taking fixed points then gives us
a non-extremal subfactor, as in \cite{Burns}.  So we obtain a continuous family of 
non-extremal subfactors with the same
principal graphs, but different indices.  

Other non-irreducible subfactors with index $(1 + \sqrt{3})^2 = 4 + 2 \sqrt{3}$ may be constructed similarly from a graph which
branches twice at each even vertex and four times at each odd vertex.  This is the third on the list from \cite{PP}.

\vspace{0.5cm}
\noindent{\it Acknowledgement.}
I would like to thank Professor Dietmar Bisch for his useful suggestions and corrections.

\end{document}